\documentclass[12pt]{elsarticle}

\usepackage{amsmath,amsthm,amssymb,dsfont,a0size,bm,indentfirst,wasysym}
\usepackage{stmaryrd,CJK,color,graphicx}
\usepackage[toc,page,title,titletoc,header]{appendix}
\usepackage{paralist}
\usepackage{epsfig}
\usepackage[colorlinks=true]{hyperref}
\hypersetup{citecolor=red}

\usepackage{remreset}

\allowdisplaybreaks

\newcommand{\bex}{\begin{eqnarray*}}
\newcommand{\eex}{\end{eqnarray*}}
\newcommand{\be}{\begin{eqnarray}}
\newcommand{\ee}{\end{eqnarray}}

\newcommand{\beex}{\begin{equation*}}
\newcommand{\eeex}{\end{equation*}}

\newcommand{\bee}{\begin{equation}}
\newcommand{\eee}{\end{equation}}
\newcommand{\bea}{\begin{aligned}}
\newcommand{\eea}{\end{aligned}}

\newcommand{\ba}{\begin{array}}
\newcommand{\ea}{\end{array}}
\newcommand{\bi}{\begin{itemize}}
\newcommand{\ei}{\end{itemize}}
\newcommand{\bn}{\begin{enumerate}}
\newcommand{\en}{\end{enumerate}}
\newcommand{\sex}[1]{\left(#1\right)}

\newcommand{\sed}[1]{\left\{#1\right\}}

\newcommand{\sen}[1]{\left\Vert#1\right\Vert}

\renewcommand{\theequation}{\arabic{equation}}

\newtheorem{thm}{\indent Theorem}

\newtheorem{defn}[thm]{\indent Definition}

\newcommand{\bp}{\begin{proof}}
\newcommand{\ep}{\end{proof}}

\newcommand{\ve}{\varepsilon}

\def\n{\nabla}

\def\lap{\triangle}

\def\calF{\mathcal{F}}

\def\calS{\mathcal{S}}

\def\calZ{\mathcal{Z}}

\def\bbu{\bm{u}}

\def\bbN{\mathbb{N}}

\def\bbR{\mathbb{R}}

\def\bbZ{\mathbb{Z}}

\def\rd{\mathrm{\,d}}

\allowdisplaybreaks
\def\ve{\varepsilon}

\makeatletter
\def\@eqnnum{{\normalfont \color{red} (\theequation)}}
\makeatother

\makeindex

\begin{document}

\begin{frontmatter}

\title{A smallness regularity criterion for the $3$D Navier-Stokes equations in the largest class}

\author[a]{Zujin Zhang\fnref{z}\corref{cor}}

\address[a]{
School of Mathematics and Computer Science,
Gannan Normal University\\
Ganzhou 341000, P.R. China}

\cortext[cor]{Corresponding author}

\fntext[z]{zhangzujin361@163.com}

\begin{abstract}
In this paper, we consider the three-dimensional Navier-Stokes equations, and show
that if the $\dot B^{-1}_{\infty,\infty}$-norm of the velocity field is sufficiently
small, then the solution is in fact classical.
\end{abstract}


\end{frontmatter}

\section{Introduction}
\label{sect:intro}
    Consider the following three-dimensional ($3$D) Navier-Stokes equations:
    \bee
    \label{NSE}
    \left\{
    \ba{l}
    \bbu_t+(\bbu\cdot\n)\bbu-\lap \bbu+\n \pi=0,\\
    \n\cdot \bbu=0,\\
    \bbu(x,0)=\bbu_0,
    \ea\right.
    \eee
    where $\bbu=(u_1(x,t),u_2(x,t),u_3(x,t))$ is the fluid velocity,
    $\pi=\pi(x,t)$ is a scalar pressure; and $\bbu_0$ is the prescribed
    initial velocity filed satisfying the compatibility condition $\n\cdot\bbu_0=0$.

    The existence of a global weak solution
    \bex
    \bbu\in L^\infty(0,T;L^2(\bbR^3))
    \cap
    L^2(0,T;H^1(\bbR^3))
    \eex
    to \eqref{NSE} has long been established by Leray
    \cite{Leray}, see also Hopf \cite{Hopf}.
    But the issue of regularity and uniqueness of $\bbu$ remains open.
    Initialed by Serrin \cite{Serrin_62,Serrin_63} and
    Prodi \cite{Prodi_59}, there have been a
    lot of literatures devoted to finding sufficient conditions
    to ensure $\bbu$ to be smooth, see, e.g., \cite{Beirao da Veiga_95,Beirao da Veiga_02,CT10_p3u3,CL01_p,
    Constantin_Fefferman_93,Escauriaza_Seregin_Sverak_03,FJN08_p,
    He_Gala_11_NARWA,Neustupa_Penel,Penel_Pokorny,Zhang_Chen_05,Zhang_arXiv_Hess,Zhang_arXiv_u3,
    ZhYLGL,Zhou_05_one component,ZP09_pu3,ZP10_u3_p3u3} and references cited therein. Noticeably,
    the following Ladyzhenskaya-Prodi-Serrin condition
    (\cite{Escauriaza_Seregin_Sverak_03,Prodi_59,Serrin_62,Serrin_63}):
    \be\label{reg_LSP}
    \bbu\in L^p(0,T;L^q(\bbR^3)),\mbox{ with }\frac{2}{p}+\frac{3}{q}=1,\quad
    3\leq q\leq \infty
    \ee
    can ensure the smoothness of the solution.

    Note that the limiting case $L^\infty(0,T;L^3(\bbR^3))$ in \eqref{reg_LSP}
    does not fall into the framework of standard
    energy method, which was proved by Escauriaza, Seregin and
    \v Sver\'ak \cite{Escauriaza_Seregin_Sverak_03} using backward uniqueness theorem.
    Due to the fact that
    \bex
    L^3(\bbR^3)\subset \dot B^{-1}_{\infty,\infty}(\bbR^3),\mbox{ but }
    L^3(\bbR^3)\neq \dot B^{-1}_{\infty,\infty}(\bbR^3),
    \eex
    we shall consider in this paper the regularity of solutions of
    \eqref{NSE} in $\dot B^{-1}_{\infty,\infty}(\bbR^3)$.
    However, we could not prove a regularity criterion as $L^\infty(0,T;\dot B^{-1}_{\infty,\infty}(\bbR^3))$,
    since the function in $\dot B^{-1}_{\infty,\infty}(\bbR^3)$ has no decay at infinity, which ensures that the
    solution is smooth outside an big ball centered at origin so that the backward uniqueness theorem can be applied.

    Before we state the precise result, let us recall the weak formulation of \eqref{NSE}.

    \begin{defn}\label{defn:wf}
    Let $\bbu_0\in L^2(\bbR^3)$ satisfying $\n\cdot\bbu_0=0$,
    $T>0$. A measurable vector-valued function $\bbu$
    defined in $[0,T]\times\bbR^3$ is said to be a weak solution to \eqref{NSE} if
    \bn
    \item $\bbu\in L^\infty(0,T;L^2(\bbR^3))\cap L^2(0,T;H^1(\bbR^3))$;
    \item $\bbu$ satisfies $\eqref{NSE}_{1,2}$ in the sense of distributions;
    \item $\bbu$ satisfies the energy inequality:
        \bex
        \sen{\bbu(t)}_{L^2}^2
        +2\int_0^t \sen{\n\bbu(s)}_{L^2}^2\rd s
        \leq \sen{\bbu_0}_{L^2},\quad\mathrm{a.e.}\quad t\in [0,T].
        \eex
    \en
    \end{defn}

    Now, our main result reads:
    \begin{thm}\label{thm:main}
    Let $\bbu_0\in L^2(\bbR^3)$ satisfying $\n\cdot\bbu_0=0$,
    $T>0$. Assume that $\bbu$ is a weak solution of \eqref{NSE} in $[0,T]$.
    If there exists an absolute constant $\ve_0>0$ such that
    \be\label{thm:main:reg}
    \sen{\bbu}_{\dot B^{-1}_{\infty,\infty}}\leq\ve_0,
    \ee
    then $\bbu$ is smooth in $(0,T)$.
    \end{thm}

    The rest of this paper is organized as follows. In section \ref{sect:Pre},
    we recall the definition of
    Besov spaces and an interpolation inequality.
    Section \ref{sect:proof} is devoted to proving  Theorem \ref{thm:main}.

\section{Preliminaries}\label{sect:Pre}

    We first introduce the Littlewood-Paley decomposition. Let $\calS(\bbR^3)$ be the Schwartz class of rapidly decreasing functions. For $f\in \calS(\bbR^3)$, its Fourier transform $\calF f=\hat f$ is defined as
    \bex
    \hat f(\xi)=\int_{\bbR^3}f(x)e^{-ix\cdot \xi}\rd x.
    \eex

    Let us choose an non-negative radial function $\varphi\in \calS(\bbR^3)$ such that
    \bex
    0\leq \hat \varphi(\xi)\leq 1,\quad \hat \varphi(\xi)=\left\{\ba{ll}
    1,&\mbox{if }|\xi|\leq 1,\\
    0,&\mbox{if }|\xi|\geq 2,
    \ea\right.
    \eex
    and let
    \bex
    \psi(x)=\varphi(x)-2^{-3}\varphi(x/2),\
    \varphi_j(x)=2^{3j}\varphi(2^jx),\
    \psi_j(x)=2^{3j}\psi(2^jx),\quad j\in\bbZ.
    \eex
    For $j\in\bbZ$, the Littlewood-Paley projection operators $S_j$ and $\lap_j$ are, respectively, defined by
    \bex
    S_jf=\varphi_j*f,\quad \lap_jf=\psi_j*f.
    \eex
    Observe that $\lap_j=S_j-S_{j-1}$. Also, it is easy to check that if $f\in L^2(\bbR^3)$, then
    \bex
    S_jf\to 0,\mbox{ as }j\to -\infty;\quad
    S_jf\to f,\mbox{ as }j\to \infty,
    \eex
    in the $L^2$ sense. By telescoping the series, we have the following Littlewood-Paley decomposition
    \bex
    f=\sum_{j=-\infty}^\infty\lap_jf,
    \eex
    for all $f\in L^2(\bbR^3)$, where the summation is in the $L^2$ sense.

    Let $s\in\bbR$; $p,q \in [1,\infty]$, the homogeneous Besov space $\dot B^s_{p,q}(\bbR^3)$ is defined by the full dyadic decomposition such as
    \bex
    \dot B^s_{p,q}
    =\sed{f\in \calZ'(\bbR^3);\ \sen{f}_{\dot B^s_{p,q}}=\sen{\sed{
    2^{js}\sen{\lap_jf}_{L^p}
    }_{j=-\infty}^\infty}_{\ell^q}<\infty},
    \eex
    where $\calZ'(\bbR^3)$ is the dual space of
    \bex
    \calZ(\bbR^3)=\sed{f\in\calS(\bbR^3);\
    D^\alpha \hat f(0)=0,\quad \forall\ \alpha\in \bbN^3}.
    \eex

    The following interpolatin inequality will be need in Section \ref{sect:proof},
    \bee\label{inter_ineq}
    \sen{f}_{L^q}
    \leq C
    \sen{f}_{\dot H^{\alpha\sex{\frac{q}{2}-1}}}^\frac{2}{q}
    \sen{f}_{\dot B^{-\alpha}_{\infty,\infty}}^{1-\frac{2}{q}},\quad
    \forall \ f\in \dot H^{\alpha\sex{\frac{q}{2}-1}}(\bbR^3)
    \cap \dot B^{-\alpha}_{\infty,\infty}(\bbR^3),
    \eee
    where $2<q<\infty$ and $\alpha>0$. See \cite{MGO_97} for the proof.

\section{Proof of Theorem \ref{thm:main}}\label{sect:proof}

    In this section, we shall prove Theorem \ref{thm:main}.

    By the classical ``weak=strong" type uniqueness theorem, we need only to derive the a priori estimate
    \bee\label{goal}
    \bbu\in L^\infty(0,T;H^1(\bbR^3))\cap L^2(0,T;H^2(\bbR^3)).
    \eee

    Multiplying $\eqref{NSE}_1$ by $-\lap \bbu$, integrating over $\bbR^3$, we obtain
    \bee\label{energy}
    \bea
    \frac{1}{2}
    \frac{\rd}{\rd t} \sen{\n\bbu}_{L^2}^2
    +\sen{\lap\bbu}_{L^2}^2
    &=\int_{\bbR^3} [(\bbu\cdot\n)\bbu]\cdot\lap\bbu\rd x\\
    &\equiv I.
    \eea
    \eee
    By H\"older inequality,
    \bex
    I\leq \sen{\bbu}_{L^6}\sen{\n\bbu}_{L^3}
        \sen{\lap\bbu}_{L^2}.
    \eex
    Invoking \eqref{inter_ineq} with $q=6$, $\alpha=1$; and $q=3$, $\alpha=2$, we may further estimate $I$ as
    \bee\label{I}
    \bea
    I&\leq C\sex{\sen{\bbu}_{\dot H^2}^\frac{1}{3}
        \sen{\bbu}_{\dot B^{-1}_{\infty,\infty}}^\frac{2}{3}}
        \sex{\sen{\n\bbu}_{\dot H^1}^\frac{2}{3}\sen{\n\bbu}_{\dot B^{-2}_{\infty,\infty}}^\frac{1}{3}}
        \sen{\lap\bbu}_{L^2}\\\
    &=C\sen{\bbu}_{\dot B^{-1}_{\infty,\infty}}\sen{\lap\bbu}_{L^2}^2.
    \eea
    \eee
    Substituting \eqref{I} into \eqref{energy}, we see
    \bex
    \frac{1}{2}
    \frac{\rd}{\rd t} \sen{\n\bbu}_{L^2}^2
    +\sex{1-C\sen{\bbu}_{\dot B^{-1}_{\infty,\infty}}}\sen{\lap\bbu}_{L^2}^2
    \leq 0.
    \eex
    Thus, if
    \bex
    \sen{\bbu}_{\dot B^{-1}_{\infty,\infty}}\leq \frac{1}{C}\equiv \ve_0,
    \eex
    we deduce that $\sen{\n\bbu}_{L^2}$ is decreasing, and thus \eqref{goal}, as desired.

    The proof of Theorem \ref{thm:main} is completed.

\end{document}